\documentclass[12pt,a4paper,oneside]{book}
\usepackage{amsmath}
\usepackage{amssymb}
\usepackage{amsthm}
\usepackage{amsfonts}
\usepackage{graphicx}
\usepackage[all]{xy}
\usepackage{verbatim}
\usepackage[left=4cm,right=2cm,top=3cm,bottom=2.5cm]{geometry}
\usepackage{caption}
\usepackage{subcaption}
\usepackage{float}
\usepackage{stackrel}
\usepackage{setspace}
\usepackage{hyperref}
\usepackage[utf8]{inputenc}
\usepackage[english]{babel}\usepackage{amsmath}
\usepackage[numbers]{natbib}
\usepackage[utf8]{inputenc}
\usepackage[english]{babel}
\usepackage{eucal}
\usepackage{float}
\usepackage{color}

\usepackage{array}
\newcolumntype{C}{@{}c@{}}
\newcommand{\bottombox}[1]{\makebox[2em][r]{#1}\hspace*{\tabcolsep}\hspace*{2em}}%
\newcommand{\innerbox}[2]{%
	\begin{tabular}[b]{c|c}
		\rule{2em}{0pt}\rule[-2ex]{0pt}{5ex} & \makebox[2em]{#2} \\\cline{2-2}
		\multicolumn{2}{r}{{#1}\hspace*{1.5\tabcolsep}\hspace*{2em}\rule[-2ex]{0pt}{5ex}}
\end{tabular}}

\begin{document}
	
	\pagenumbering{roman}

	\pagestyle{headings}

	\begin{center}
		
		\vspace{0.5in}
		{\bf{Application of KKT to Transportation Problem with Volume Discount on Shipping Cost.}\\[5mm] Haruna Issaka\textsuperscript{1}, Ahmed mubarack\textsuperscript{2}, Osman Shaibu\textsuperscript{3}}\\PhD Computational Mathematics.\\2. PhD Student, University of Nottingham, United Kingdom.\\\underline{3. Department of Mathematics, University of Health and Allied Sciences, Ghana.}
	\end{center}

	\setcounter{page}{1}
	\pagenumbering{arabic}
	\subsection*{Abstract} 
	Nonlinear programming problems are useful in designing and assigning work schedule and also in transporting goods and services from known sources to specified destinations. The objective function could be linear or nonlinear depending on the mode or rout of transportation. In supplying goods and services during emergency cases, the objective cost function is always assumed to be nonlinear. 
	\\
	The purpose of this research is to study the nature of the cost function of the transportation problem and provide solution algorithms to solve such nonlinear cost functions. The cost per unit commodity transported may not be fixed since volume discounts are sometimes available for large shipments of goods and services. This will make the cost function either piecewise linear, separable concave or convex.  
	The approach to solving such problems as they arise is to apply existing general nonlinear programming algorithms and when necessary make modifications to fit the special structure of the problem. We describe theoretically , algorithms for the nonlinear cost function and state conditions under which stability is achieved with both sufficient and necessary cases.	
	
\noindent	\underline{\textit{keywords:} Convex function, shipping cost, transportation problem,\\ optimization problem.}

	\section*{1. Introduction} A number of factors needs to be taken into account when considering transportation problems. For instance, port selection, rout of transport, port to port carrier selection as well as distribution-related cases. Freight companies dealing with large volume shipments, encounter various problems during the movement of good and services due to the cost involved. Volume discounts are meant to target minimizing shipping cost of transport \cite{issaka2010transportation}. A very important area of application is the management and efficient use of scarce resources to increase productivity \cite{adby1974introduction},\cite{mccormick1983nonlinear}. 
	\\
	\cite{shetty1959solution} formulated an algorithm to solve transportation problems taking nonlinear cost function. 
 In his paper, solution to the generalized transportation problem taking nonlinear costs is given. The solution procedure used is an iterative method and a feasible solution is obtained at each stage. The value of the criterion function is improved from one stage to another.
 \\
 Bhudury et al \cite{arcelus1996buyer}  gave a detailed game-theory analysis of the buyer-vendor coordination problem embedded in the price-discount inventory model. Pure and mixed, cooperative and non-cooperative strategies were developed. Highlights of the paper include the full characterization of the Pareto optimal set, the determination of profit-sharing mechanisms for the cooperative case and the derivation of a set of parameter-specific non-cooperative mixed strategies. Numerical example was presented to illustrate the main features of the model.
 \\
 \cite{goyal1988joint} also developed a joint economic‐lot‐size model for the case where a vendor produces according to order for a purchaser on a lot‐for‐lot basis under deterministic conditions. The assumption of lot‐for‐lot bases is restrictive in nature. In his paper, general joint economic‐lot‐size model was suggested and it was shown to provide a lower or equal joint total relevant cost as compared to the model of Banerjee.
 \\
 In their paper, \cite{sharp1970decomposition},\cite{abramovich2004refining} considered a firm with several plants supplying several markets with known demands. The cost of production was considered to be nonlinear while transportation cost between any plant and market pair are linear. The Kuhn-Tucker conditions and the dual to the transportation model were used to derive the optimal conditions for the problem. The conditions were shown to be both sufficient and necessary if the production cost are convex at each plant else it is only necessary. An algorithm was also developed for reaching an optimal solution to the production-transportation problem for the convex case. 
 \\
 Tapia et al \cite{tapia1994extension} extended the Karush–Kuhn–Tucker approach to Lagrange multiplier of an infinite programming formulation. The main result generalized the usual first-order necessity conditions to address problems in which the domain of the objective function is Hilbert space and the number of constraints is arbitrary. 
 Under production economies of scale, the trade‐off between production and transportation costs will create a tendency towards a decentralized decision‐making process \cite{youssef1996iterative}.
 \\
 The algorithms of this paper belong to the direct-search or implicit-enumeration type algorithms was proposed by  \cite{spielberg1969algorithms}. In the paper, a general plan of procedure is expected to be equally valid for the capacitated plant-location problem and also for transshipment problems with fixed charges. After considerable computational experience accumulated and discussed at some length, the authors suggested an additional work on the construction of adaptive programs that matches algorithm to data structure.
 \\
  Williams, \cite{williams1962treatment} applied the decomposition principle of \cite{dantzig1959truck} to the solution of the \cite{hitchcock1941distribution} transportation problem and to several generalizations of it. Among the generalizations are $(i)$ the transportation problem with source availability subject to general linear constraints, and $(ii)$ the case in which the costs are piecewise linear convex functions.

	\subsection*{2. Non-Linear programming Problem}
	In order to study nonlinear problems, we observe the following definitions:
	\\
	$i.$ \textbf{Polyhedral set}: A set $P$ in an $n$ dimensional vector space $T^n$ is a polyhedral set if it is the intersection of a finite number of closed- half	spaces, i.e. $P = {x : p^t_ix}\leq \alpha _i \quad \forall i=1...n$, where $p_i$ is a non zero vector in $T^n$ and $\alpha_i$ is scalar.
	\\
	$ii$.  \textbf{Convex set}: A subset $S$ of $R^n$ is convex if for any $x_1, x_2 \in S$, the line segment $[x_1, x_2]$ is contained in $S$. In other words the set $S$ is convex if $\lambda x_1+(1-\lambda)x_2\in S$ \quad $\forall \quad 0\leq \lambda \leq 1$.
	\subsection*{3. Karish Kuhn Tucker conditions for Optimality}
	Consider the non-linear transportation 
	\\
	Min $C_i(x_i)$ subject to the constraint $A x_i=b_i$\quad $\forall x_i>0$ and $i=1,2,...$, where $x=(x_1,x_2,...x_n)^T$,
	$b=\begin{pmatrix}
	s_{1}\\s_{2}\\..\\{s_n}\\{d_1}\\..\\{d_n}
	\end{pmatrix}$, $A=\begin{pmatrix}
	1 &1&.&.&.&1 \\ .&.&.&.&.&.&\\1&1&.&.&.&1
	\end{pmatrix}$.
	\\
	$C_i$ is the cost associated with each $x_i$.
	\\
	The solution tableau is setup as below
	\\[5mm]
	    \begin{tabular}{|c|C|C|C|r|}\hline
	    	& $x_{11}$                 & $x_{12}$                 & $x_{ij}$                 & Supply  \\\hline
	    	1      & \innerbox{$x_{1,1}$}{$C_{1,1}$}   & \innerbox{$x_{1,2}$}{$C_{1,2}$}   & \innerbox{$x_{1,j}$}{$C_{1,j}$} &  $s_1$   \\\hline
	    	2      & \innerbox{$x_{2,1}$}{$C_{2,1}$} & \innerbox{$x_{2,2}$}{$C_{2,2}$}   & \innerbox{$x_{2,j}$}{$C_{2,j}$}    &   $s_2$   \\\hline
	    	..      & \innerbox{..}{..} & \innerbox{..}{..}    & \innerbox{..}{..}    &   ..  \\\hline
	    	n      & \innerbox{$x_{1n}$}{$C_{1n}$}   & \innerbox{$x_{2,n}$}{$C_{2,n}$} & \innerbox{$x_{i,n}$}{$C_{i,n}$} &   $s_n$  \\\hline
	    	n+1  & \innerbox{$x_{1,n+1}$}{$C_{1,n+1}$} & \innerbox{$x_{2,n+1}$}{$C_{2,n+1}$}    & \innerbox{$x_{i,n+1}$}{$C_{i,n+1}$}    &   $s_{n+1}$   \\\hline
	    	Demand & \bottombox{$d_1$}   & \bottombox{$d_2$}   & \bottombox{$d_n$}   & $\sum_{i=1}^{n}s_i=\sum_{j=1}^{n}d_j$   \\\hline
	    \end{tabular}
	\\[5mm]
	where $C_{ij}$ is the cost of transporting quantity $x_{ij}$ from source $s_i$ to destination $d_j$.
	\\
	The Lagrangian function for the system is $Z(x,\lambda,\omega)=Cx+\omega (b-Ax)-\lambda x$ where $\lambda$ and $\omega$ are Lagrangian multipliers. By the KKT condition, the basic solution $\bar{x}$ is feasible if \\
	\begin{equation}
	\begin{cases}
	\label{eqn1}
	\triangledown Z=\triangledown C(\bar{x})-\omega^TA-\lambda=0,\\ \lambda \bar{x}=0,\\ \lambda \geq 0, \\ \bar{x}\geq 0.
	\end{cases}
	\end{equation}
	
	from \ref{eqn1}, we have the following
	\\
	\begin{equation}
	\label{eqn2}
	\begin{cases}
	\frac{\partial Z}{\partial x_{ij}}=\frac{\partial C(\bar{x})}{\partial x_{ij}}-(u,v)(e_i,e_{n+j})-\lambda_k=0,
	\\
	\lambda_{ij}x_{ij}=0\\
	x_{ij}\geq 0\\
	\lambda_k\geq 0.
	\end{cases}
	\end{equation}
	
\noindent	From \ref{eqn2} we have \\
	\begin{equation}
	\label{eqn3}
	\frac{\partial Z}{\partial x_{ij}}=\frac{\partial C(\bar{x})}{\partial x_{ij}}-(u_i,v_j) \geq 0
	\end{equation}
	and
	\begin{equation}
	\label{eqn4}
	x_{ij}\frac{\partial Z}{\partial x_{ij}}=\frac{\partial C(\bar{x})}{\partial x_{ij}}-(u_i,v_j) \geq 0 \quad \forall x_{ij}\geq 0.
	\end{equation}
	 
	\noindent Hence, when ever conditions $3$ and $4$ are satisfied, the system attains its optimal value. 
	 
	\subsection*{4. General Solution to the Non-linear Transportation Problem}
	The solution process is detailed below
	\begin{enumerate}
		\item Find an initial basic feasible solution $\bar{x}$ using any of the known methods. i.e. North-west corner rule, Vogel's approximation, Row minima etc.
		\item If $\bar{x}$ satisfies \ref{eqn3} and \ref{eqn4}, then $\bar{x}$ is KKT point and therefore stop
		\item If $\bar{x}$ does not satisfies \ref{eqn3} and \ref{eqn4}, then move to find a new basic solution by improving the cost function and repeat the process.
	\end{enumerate}
	
	\subsection*{5. Non-linear Transportation Problem with Concave Cost Function}
	During emergency cases, goods and services are shipped through unapproved routs to the affected destination. In such cases, the cost function assumes a nonlinear structure. The cost function is concave in nature with a concave structure. Volume discounts are sometimes available for shipments in large volumes. When volume discounts are available, the cost per unit shipped decreases with increasing volume shipped. The volume may be directly related to the unit shipped or have the same rate for equal amounts shipped.\\
	If the discount is directly related to the unit commodity shipped, the resulting cost function will be continuous with continuous first partial derivative.\\ The associated nonlinear problem formulation is
	\\
	Min$f(x)=C_{ij}x_{ij}$ subject to the constraints \\ $\sum_{i=1}^{n} x_{ij}=s_i$ and $\sum_{j=1}^{n} x_{ij}=d_j$ where the cost function $C_{ij} : R^{mn}\rightarrow R$. 
	\\[2mm]
	\textbf{Theorem 1}: Let $f$ be concave and continuous function and $P$ a non-empty compact polyhedral set. Then the optimal solution to the problem $min f(x), \quad x\in P$ exist and can be found at an extreme point of $P$.\\
	\textbf{Proof}: Let $E=(x_1, x_2,...,x_n)$ be an extreme point of $P$ and $x_k\in E$ such that\\ $f(x_k)=min \{ f(x_i)\mid i=1,2,...n \}$. Since $P$ is compact and $f$ is continuous, $f$ attains it's minimum in $P$ at $x_k$. If $\bar{x}$ is an extreme point of $P$, the solution is optimal and $\bar{x}=\sum_{i=1}^{n}\lambda_ix_i$, $\sum_{i=1}^{n}\lambda_i=1$, $\lambda_i >0$ for all extreme points $(x_1, x_2, ...x_n)$. \\ By the concavity of $f$, it follows that\\ $f(\bar{x})=f(\sum_{i=1}^{n}\lambda_ix_i)\geq \sum_{i=1}^{n}\lambda_if(x_i)\geq f(x_k)\sum_{i=1}^{n}\lambda_i$.
	This implies $f(\bar{x})\geq f(x_k)$ since $f(x_k)\leq f(x_i)$ and $\sum_{i=1}^{n}\lambda_i=1$. Therefore $\bar{x}$ is minimum and so $f(\bar{x})\leq f(x_k)$. Hence $f(\bar{x})=f(x_k)$
	
	\subsection*{6. Solution to Concave Cost Function}
	Theorem $1$ allows us to solve the concave problem as follows:\\ We consider only the extreme points of $P$ to minimize the cost $C_{ij}$.\\ Let $\bar{x}$ be the basic solution in the current iteration. We decompose $\bar{x}$ into $(x_B, x_N)$ where $x_B$ and $x_N$ are respectively the basic and non-basic variables.
	\\
	 For $x_B>0$, the complementary slackness condition is $m+n-1$, where $m$ and $n$ are the demand and supply capacities.
	 \\
	Next we determine the values of $u_i$ and $v_j$ from\\ $\frac{\partial C_{ij}x_{ij}}{\partial x_{Bij}}-(u_i+v_j)=0$.
	\\
    We determine $u_i$ and $v_j$ for all $i,j=1,...n$ by assigning $u_1$ the value zero and solving for the remaining $u_i's$ and $v_j's$. for each non-basic variable $x_{Nij}$, we calculate from \ref{eqn3}, $\frac{\partial Z}{\partial x_ij}$.
    \\
     At the extreme point of $P$, the $x_{ij}'s$ are zeros and hence the complementary slackness condition is satisfied.
     \\
	On the other hand, if 
	\\
	$\frac{\partial Z}{\partial x_ij}-(u_i+v_j)<0$ for all $ij=1,...,n$, we move to search for a better basic solution by determining the leaving and new entering basic variables.
	
	\subsection*{7. Convex Cost Function}
	\textbf{Definition:} A line segment joining the points $x$ and  $y$ in $R^n$ is the set ${x,y}$ = ${x\in R^n:x=\lambda x + (1 – \lambda)y}$, $\forall$ $0 \leq \lambda \leq 1$. A point on the line segment for which $0 < \lambda <1$, is called an interior point of the line segment.\\
	A subset $S$ of $R^n$ is said to be convex if for any two elements $x$, $y$ in $S$,  the line segment $[x,y]$ is contained in $S$ . Thus  $x$ and $y$ in $S$ imply ${x + (1 – \lambda)y}$ $\in S$ for $0 \leq \lambda \leq 1$ if $S$ is convex.
	
	\subsection*{8. Convex Optimization Problem:}
	A convex optimization problem is an optimization technique that consist of minimizing a convex cost function over a convex set of constraints.
	\\
	 Mathematically we have, 
	\textit{min} $f(x)$ subject to $x\in C$, where $C$ is a convex set and $f$ a convex function in $C$. \\
	In particular, \textit{min} $f(x)$ subject to $g_i\leq 0$, $i=1,2,..n$, where $f$ is a convex function and $g$ the associated constraint.\\ In line with the above, we state the following theorem.
	\\[4mm]
	$\textbf{Theorem 2:}$
	\\
	 Let $f:C\mapsto R$ be convex function on the set $C$. If $\bar{x}$ is a local minimum of $f$ over $C$, i.e. $\bar{x}$, then $\bar{x}$ is a global minimum of $f$ over $C$. 
	\\
\noindent	\textbf{Proof:}\\
	If $\bar{x}$ is a local minimum of $f$ over $C$, then there exist $r>0$ such that $f(x)\geq f(\bar{x})$ for all $x\in C$. Let $y\in C$ such that $y\neq \bar{x}$. It suffices to show that $f(y)\geq f(\bar{x})$. Suppose $\lambda \in (0,1]$ is such that $\bar{x}+\lambda(y-\bar{x})\in [\bar{x}, r]$. It follows from Jensen’s inequality that $f(\bar{x})\leq f(\bar{x}+\lambda (y-\bar{x}))\leq (1-\lambda)f(\bar{x})-\lambda f(y)$.
	\\
	$\implies \lambda f(\bar{x})\leq\lambda f(y)$ and hence $f(\bar{x})\leq f(y)$. 
	
	\subsection*{9. Convex Transportation Problem:}
	During emergency situations (natural disasters like flooding) where logistics are required to be moved immediately to site, the mode and method of transportation is sometimes not regular. In developing a cost function for such a case, the objective function assumes a nonlinear form.
	This case arise when the objective function is composed of not only the unit transportation cost but also of production cost related to each commodity. Or in the case when the distance from each source to each destination is not fixed, the objective function nonlinear; concave or convex cost function. 
	The associated problem is formulated as:
	
\noindent	\textit{Minimize} $C(x)$\\ subject to $Ax=b$ $x\geq0$\\ where $C(x)$ is convex, continuous and has continuous first order partial derivatives.
	
	\subsection*{10. Solution to the Convex Cost function}
	For a convex cost function, the minimum point may not be attained at an extreme point of the function but a solution may be reached at the boundary of the feasible solution set. This case may arise when a nonbasic variable has a positive allocation while none of the basis is driven to zero.\\
	We use Zangwill's \cite{abramovich2004refining} convex simplex algorithm to seek for a local optimal solution by partitioning the solution variable $x$ into $\{x_B, x_N\}$ where $x_B$ is the $n+m-1$ component of the basic variable and $x_N$ is the $nm-(nm-1)$ component vector of the nonbasic variable.\\
	Suppose we have an initial basic solution $\bar{x}$, we improve the solution by the KKT conditions studied earlier, until each cell $x_{ij}$ satisfies sufficiently the following conditions:
	\begin{enumerate}
		\item $x_{ij}\dfrac{\partial f(\bar{x})}{\partial x_{ij}}-(u_i+v_j)=0$,
		\item $\dfrac{\partial f(\bar{x})}{\partial x_{ij}}-(u_i+v_j)\geq0$
	\end{enumerate}    
	
\noindent	Since each $x_{B_{ij}}>0$, condition $2$, the slackness condition, becomes\\
	
	 $\dfrac{\partial f(\bar{x})}{\partial x_{B_{ij}}}-(u_i+v_j)=0$ where $x_{B_{ij}}$ is a basic variable in the $ij$ cell.
	\\
	But for a nonbasic cell $x_N$ at a feasible iteration point, any of the following may occur:
	\\
	\begin{enumerate}
		\item $\dfrac{\partial f(\bar{x})}{\partial x_{ij}}-(u_i+v_j)>0$ \quad
		 and \quad $x_{ij}\dfrac{\partial f(\bar{x})}{\partial x_{ij}}-(u_i+v_j)>0$
		\item $\dfrac{\partial f(\bar{x})}{\partial x_{ij}}-(u_i+v_j)<0$ \quad
		 and \quad $x_{ij}\dfrac{\partial f(\bar{x})}{\partial x_{ij}}-(u_i+v_j)<0$
		\item $\dfrac{\partial f(\bar{x})}{\partial x_{ij}}-(u_i+v_j)<0$ \quad
		and \quad $x_{ij}\dfrac{\partial f(\bar{x})}{\partial x_{ij}}-(u_i+v_j)=0$
	\end{enumerate}
	
\noindent	If $x$ satisfies any of the three conditions, we improve the solution as follows:
	\\
	We compute $\frac{\partial z}{\partial x_{rl}}$, the minimum of $\dfrac{\partial f(\bar{x})}{\partial x_{ij}}-(u_i+v_j)$ and 
	\\
	 $x_{st}\frac{\partial z}{\partial x_{st}}$ the maximum of $x_{ij}\dfrac{\partial f(\bar{x})}{\partial x_{ij}}-(u_i+v_j)$
	 \\
	 At this point we only improve a nonbasic variable $x_{ij}$ with a positive value and a positive partial derivative.\\ We select the leaving variable under each of the three conditions stated above as follows:
	 \\
	 \textbf{Case 1}
	 \\
	 If $\frac{\partial z}{\partial x_{rl}}\geq0$	 and $x_{st}\frac{\partial z}{\partial x_{st}}>0$, we decrease $x_st$ by the value $\theta$ using the transportation table as in the linear and concave cases.
	 Let $y^k = (y_{11}^k, y_{12}^k,...,y_{mn}^k)$ be the value of $x^k = (x_{11}^k,..., x_{mn}^k)$ after making the necessary adjustment by adding and subtracting $\theta$ in the loop containing $x_{st}$ so that all the constraints are satisfied.
	 By doing so, either $x_{st}$ itself or a basic variable say $x_{B_{st}}$ will be driven to zero.
	 Now $y^k$ may not be the next iteration point. Since the function is convex, a better point could be found before reaching $y^k$. To check this, we solve 
	 
	\begin{equation}
	\label{eqn5}
	 f(x^{k+1}) = min f\lambda x^k + (1–\lambda)y^k) : 0 \leq \lambda \leq 1
	\end{equation}
	 
\noindent	 and obtain $x^{k+1}=\lambda x^k + (1 –\lambda)y^k$ where $\lambda$ is the optimal solution of \ref{eqn5}.
	 Before the next iteration, If  $x^{k+1} = y^k$ and a basic variable goes to zero, we change the basis.
	 If $x^{k+1} \neq y^k$ or  if $x^{k+1} = y^k$ and $x_st$ is driven to zero, we don’t change the basic by substituting the leaving basic variable by $x_st$.
	 \\
	 \textbf{Case 2}
	 \\
	 If $\frac{\partial z}{\partial x_{rl}}\leq 0$ and $x_{st}\frac{\partial z}{\partial x_{st}}\leq 0$
	 In this case the value of $x_rl$ is increased by $\theta$ and we find $y^k$, where $\theta$ and $y^k$ are defined as in case $1$.
	 As we increase the value of $x_rl$ one of the basic variables will be driven to zero.	 
	 But now after solving for $x^k+1$ from \ref{eqn5}, before going to the next iteration, we will have the following possibilities.
	 If $x^{k+1} = y^k$ we change the former basis, substitute $x_B$ by $x_{rl}$.
	 If $x^{k+1} \neq y^k$ we do not change the basis.
	 All the basic variables outside of the loop will remain unchanged.
	 Case $3$ 
	 If	$\frac{\partial z}{\partial x_{rl}}<0$ and $x_{st}\frac{\partial z}{\partial x_{st}}>0$,	we either decrease $x_{st}$ as in the case $1$ or increase $x_{rl}$ as in case $2$.
	 
	 \subsection*{11. The Transportation Convex Simplex Algorithm}
	 Algorithm for solving the convex transportation problem is developed in the following three steps.
	 \\
	 \textbf{Initialization}
	 \\
	 Final the initial basic feasible solution
	 Iteration
	 \\
	 \textbf{Step 1}: We determine all $u_i$ and $v_j$ from $\frac{\partial f(x^k)}{\partial x_{B_{rl}}}-u_i-v_j=0$
	 for each basic cell.
	 \textbf{Step 2}:  For each non basic cell, calculate
	 $$\frac{\partial z}{\partial x_{rl}}=min\{\frac{\partial f(x^k)}{\partial x_{ij}}-u_i-v_j\}$$
	 $$x_{st}\frac{\partial z}{\partial x_{rl}}=max\{x_{ij}\frac{\partial f(x^k)}{\partial x_{ij}}-u_i-v_j\}$$
	 \\
	 If
	 $\frac{\partial z}{\partial x_{rl}}\geq 0$ and $x_{st}\frac{\partial z}{\partial x_{st}}=0$
	we stop otherwise we move to step $3$.
	 \\
	 \textbf{Step 3}:
	 Determine the non basic variable to change.
	 Decrease $x_st$ according to case $1$ if
	 $\frac{\partial z}{\partial x_{rl}}< 0$ and $x_{st}\frac{\partial z}{\partial x_{st}}>0$.
	 Increase $x_rl$ according to case $2$ if 
	 $\frac{\partial z}{\partial x_{rl}}< 0$ and $x_{st}\frac{\partial z}{\partial x_{st}}\leq0$
	 Either increase $x_rl$ or decrease $x_st$ if
	 $\frac{\partial z}{\partial x_{rl}}< 0$ and $x_{st}\frac{\partial z}{\partial x_{st}}<0$
	 Step $4$:   
	 Find the values of $y^k$, by means of $\theta$, and $x^{k+1}$, from \ref{eqn5}.
	 If  $y^k = x^{k+1}$ and a basic variable is driven to zero, change the basic otherwise do not change the basis.
	 $x^k = x^{k+1}$
	 go to step $1$.              
	 
	 \subsection*{Conclusion}
	 In conclusion, we have established, under considerable conditions, solution method for the various problems that may arise in setting up sets of equations defining a given transportation problem with its associated constraints. Unlike the traditional methods of setting up the objective function and its associated constraints which is always linear, shipments during emergency cases where supply is scheduled to reach particular destinations within a specified time period, comes with its own accompanying problems. We found out that during such cases, the objective function is nonlinear. Solution algorithms were therefore developed taking into account existing algorithms with reasonable modifications to accommodate the nonlinear objective function. These algorithms can be applied to real industry data to minimize cost of shipment.       
	
	
	\subsection*{13. Conflict of Interest.}
	The authors wish to categorically state that there is no conflict of interest what soever regarding the publication of this research work. 
	
	\renewcommand{\bibname}{References}
	\bibliographystyle{abbrvnat}
	\bibliography{references}
	\addcontentsline{toc}{chapter}{References}
	\nocite{*}
\end{document}